\documentstyle[epsf,12pt]{article}
\def\Bbb R{{\rm \bf R}}
\def\proclaim#1{\vskip2mm{\bf #1}\em}
\def\endproclaim{\em \vskip2mm}
\def\tag#1{\eqno(#1)}
\def\gathered{\begin{array}{c}}
\def\endgathered{\end{array}}
\def\text{\mbox}

\begin{document}

\title {Travel Time and Heat Equation.\\
One space dimensional case}
\author{Masaru IKEHATA\\
Department of Mathematics,
Faculty of Engineering\\
Gunma University, Kiryu 376-8515, JAPAN}
\date{10 Oct 2006}
\maketitle
\begin{abstract}
The extraction problem of information about the location and shape of the cavity from a
single set of the temperature and heat flux on the
boundary of the conductor and finite time interval is a typical and 
important inverse problem.  Its one space dimensional version is considered.  
It is shown that the enclosure method developed by the author for elliptic equations
yields the extraction formula of a quantity which can be interpreted as the {\it travel
time} of a {\it virtual} signal with an arbitrary fixed propagation speed that
starts at the known boundary and the initial time, reflects at another 
unknown boundary and returns to the original boundary.

\noindent
AMS: 35R30, 80A23

\end{abstract}

\section{Introduction}

\noindent In this paper we consider one space dimensional version of
a typical and important inverse problem for the heat equation.
Let $\Omega$ be a bounded domain of $\Bbb R^n (n=1,2,3)$ with smooth boundary.
Let $D$ be an open subset with Lipschitz boundary of $\Omega$ such that
$\overline D\subset\Omega$ and $\Omega\setminus\overline D$ is connected.
Let $T$ be an arbitrary positive number. Let $u=u(x,t)$ be an {\it arbitrary non constant} solution of
the heat conduction problem:
$$\begin{array}{c}
\displaystyle
u_t=\triangle u\,\,\text{in}\,(\Omega\setminus\overline D)\times]0,\,T[,\\
\\
\displaystyle
\frac{\partial u}{\partial\nu}=0\,\,\text{in}\,\partial D\times\,]0,\,T[,\\
\\
\displaystyle
u(x,0)=0\,\,\text{in}\,\Omega\setminus\overline D.
\end{array}
$$

\noindent
Here $\nu$ denotes the unit outward normal vector field on $\partial(\Omega\setminus\overline D)$.
The boundary condition for $\partial u/\partial\nu$ on $\partial D$ means that $D$ is perfectly
insulated.  Note that this is the simplest model of the heat conduction in a conductive body
having a hole or cavity $D$.

The problem to be solved is

{\bf\noindent Inverse Problem.}
Extract information about the location and shape of $D$ from the {\it single set}
of the {\it observation data} $u(x,t)$ and $\partial u/\partial\nu(x,t)$ on $\partial\Omega\times
\,]0,\,T[$.  Note that $T$ is {\it fixed}.

\noindent
The problem comes from the thermal imaging and the solution method may have
applications to, for example, the detection of corrosion.
There are extensive studies of uniqueness and stability issue of Inverse Problem.
In particular, it is known that the observation data uniquely determine general $D$ itself
under a suitable condition on the heat flux on $\partial\Omega$.
See Bryan-Caudill \cite{BC}, Canuto-Rosset-Vessella \cite{CRV} and
references therein for more information about the issue.

In this paper we are especially interested in seeking an {\it
analytical formula} for the purpose and start with {\it one space
dimensional version} of the problem.  In one space dimensional
case there is a way of calculating the so-called {\it response
operator} from a general single set of the observation data.
Moreover in Avdonin-Belishev-Rozhkov
\cite {ABR} it is shown that from the response
operator one can extract the spectral data.  Thus in one space
dimensional case the problem can be reduced to the {\it inverse spectral
problem} which has been studied well.  However, the reduction is
based on the {\it boundary controllability} for the heat equation and
not an easy way. To our best knowledge there is no attempt for
finding the {\it extraction formula} of information about the
location and shape of $D$ from the single set of the observation
data {\it without reducing} to other inverse problems. 
In this paper, we present
a direct approach using the idea of the {\it enclosure method} for
elliptic equations introduced by the author \cite{I1}.

So what is the enclosure method?
Consider a non constant solution of the elliptic problem:
$$\begin{array}{c}
\displaystyle
\triangle u=0\,\,\text{in}\,\Omega\setminus\overline D,\\
\\
\displaystyle
\frac{\partial u}{\partial\nu}=0\,\,\text{on}\,\partial D.
\end{array}
$$
In \cite{I1} we considered the problem of extracting information
about the location and shape of $D$ in two dimensions from the
observation data that is single set of Cauchy data of $u$ on
$\partial\Omega$.   Assuming that $D$ is given by the inside of a
polygon with an additional condition on the diameter, we
established an extraction formula of the {\it convex hull} of $D$
from the data.
In the method a special exponential solution of the Laplace
equation played the central role. The solution takes the form
$e^{-\tau\,s}e^{\tau x\cdot(\omega+i\,\omega^{\perp})}$
where $\tau(>0)$ and $s$ are parameters; both $\omega$ and $\omega^{\perp
}$ are unit vectors
and satisfy $\omega\cdot\omega^{\perp}=0$.
The solution divides the space into two half planes which
have a line $\{x\,\vert\,x\cdot\omega=s\}$ as the common boundary. In one part
$\{x\,\vert\,x\cdot\omega>s\}$ the solution is growing as
$\tau\longrightarrow\infty$ and in another part $\{x\,\vert\,x\cdot\omega<s\}$ decaying.
We define a
function $I_{\omega,\omega^{\perp}}(\tau,s)$ of the independent
variable $\tau$ with parameter $s$ which is called the
{\it indicator function} and can be calculated from the
observation data:
$$\displaystyle
I_{\omega,\omega^{\perp}}(\tau,s)
=e^{-\tau \,s}\int_{\partial\Omega}\{-\frac{\partial}{\partial\nu}e^{\tau x\cdot(\omega+i\omega^{\perp})
}u+\frac{\partial u}{\partial\nu}
e^{\tau x\cdot(\omega+i\omega^{\perp})
}\}ds.
$$
In \cite{I1} it is clarified that the asymptotic behaviour
of the indicator function as $\tau\longrightarrow\infty$ depends
on the position of half plane $x\cdot\omega>s$ relative to $D$.
Recall the support function $h_D(\omega)=\sup_{x\in\,D}x\cdot\omega$.
We say that $\omega$ is {\it regular} if the set $\{x\,\vert\,x\cdot\omega=h_D(\omega)\}\cap\partial D$
consists of only one point.
Then we established: for regular $\omega$
$I_{\omega,\omega^{\perp}}(\tau,s)\vert_{s=h_D(\omega)}$
is {\it truly} algebraic decaying as $\tau\longrightarrow\infty$.
This means that: there exist positive constants $A$ and $\mu$ such that
$$\displaystyle
\lim_{\tau\longrightarrow\infty}\tau^{\mu}\vert I_{\omega,\omega^{\perp}}(\tau,s)\vert_{s=h_D(\omega)}\vert
=A.
\tag {1.1}
$$
This fact is the core of the enclosure method.
Since we have the trivial identity
$$\displaystyle
I_{\omega,\omega^{\perp}}(\tau,t)
=e^{-\tau(t-s)}
I_{\omega,\omega^{\perp}}(\tau,s),
$$
from (1.1) one could conclude that:
if $s>h_D(\omega)$, then the indicator function is
decaying
exponentially; if $s=h_D(\omega)$, then the indicator function is decaying truly algebraically;
if $s<h_D(\omega)$, then the indicator function is growing exponentially.
Moreover by taking the logarithm of both sides of (1.1), we obtained also the {\it one line} formula
$$\displaystyle
\lim_{\tau\longrightarrow\infty}\frac{\log\vert I_{\omega,\omega^{\perp}}(\tau,0)\vert}{\tau}
=h_D(\omega).
\tag {1.2}
$$
By integration by parts (1.1) is equivalent to the statement:
the integral
$$
\displaystyle
e^{-\tau h_D(\omega)}\int_{\partial D}u\frac{\partial}{\partial\nu}e^{\tau x\cdot(\omega+i\omega^{\perp})}ds
\tag {1.3}
$$
is truly algebraic decaying as $\tau\longrightarrow\infty$.  The key point for this fact is,
roughly speaking, there does not exist a harmonic extension in a neighbourhood of
the point $x_0$ with $x_0\cdot\omega=h_D(\omega)$.

\noindent
This is the essence of the idea of the enclosure method.
However this is the case when the governing equation is
elliptic \cite{I2,I3,IO}. How can one apply the enclosure method
to non elliptic case?  This is another motivation of the study and
here we give an answer to this question.

Now state our one space dimensional problem.
Let $a>0$.
Let $u=u(x,t)$ be an arbitrary solution of
the problem:
$$\begin{array}{c}
\displaystyle
u_t=u_{xx}\,\,\text{in}\,]0,\,a[\times]0,\,T[,\\
\\
\displaystyle
u_x(a,t)=0\,\,\text{for}\,t\in\,]0,\,T[,\\
\\
\displaystyle
u(x,0)=0\,\,\text{in}\,]0,\,a[.
\end{array}
\tag {1.4}
$$

\noindent
Then the problem is: extract $a$ from $u(0,t)$ and $u_x(0,t)$ for $0<t<T$.

\noindent
Let $c$ be an arbitrary positive number. Let
$$\displaystyle
v(x,t)=e^{-z^2t}\,e^{xz}
\tag {1.5}
$$
where $\tau$ satisfies $\tau>c^{-2}$ and
$$\displaystyle
z=-c\tau\left(1+i\sqrt{1-\frac{1}{c^2\tau}}\,\right).
$$

\noindent
The function $v$ is a complex valued function and satisfies the backward heat equation $v_t+\triangle v=0$.
Moreover $e^{\tau s}v$ has the special character

\noindent
$\bullet$ if $s<cx+t$, then $\lim_{\tau\longrightarrow\infty}e^{\tau s}\vert v(x,t)\vert=0$

\noindent
$\bullet$ if $s>cx+t$, then $\lim_{\tau\longrightarrow\infty}e^{\tau s}\vert v(x,t)\vert=\infty$.

\noindent
The half plane $\{(x,t)\,\vert\,s>cx+t\}$ in the space time plays the same role as the half plane
$\{x\,\vert\,x\cdot\omega>s\}$ for elliptic case.  Changing $c$ means changing normal vector
of the line $cx+t=s$ and corresponds to changing $\omega$.

{\bf\noindent Definition 1.1.}
Given $c>0$, $s\in\Bbb R$,
define the {\it indicator function}
$I_{c}(\tau;s)$ by the formula
$$\displaystyle
I_{c}(\tau;s)
=e^{\tau s}
\int_0^T\left(-v_x(0,t)u(0,t)+u_x(0,t)v(0,t)\right)dt,\,\,\tau>c^{-2}
$$
where $u$ satisfies (1.4) and $v$ is the function given by (1.5).

Our main result is the following extraction formula.

\proclaim{\noindent Theorem 1.1.} Assume that we know a positive
number $M$ such that $M\ge 2a$. Let $T$ and $T'$ satisfy $0<T'\le
T$. Let $c$ be an arbitrary positive number satisfying
$Mc<\min\,\{T,2T'\}$. Assume that $u_x(0,t)$ coincides with a
polynomial of $t$ on the interval $]0,\,T'[$ one of whose
coefficients are not $0$. The formula
$$\displaystyle
\lim_{\tau\longrightarrow\infty}
\frac{\displaystyle\log\vert I_{c}(\tau;0)\vert}
{\tau}
=-2ca
\tag {1.6}
$$
and the following statements are true:

if $s\le 2ca$, then
$\lim_{\tau\longrightarrow\infty}\vert I_{c}(\tau;s)\vert=0$;

if $s>2ca$, then
$\lim_{\tau\longrightarrow\infty}\vert I_{c}(\tau;s)\vert=\infty$.

\endproclaim

\noindent
This is an unexpected result.
Under the condition on $T$
integration by parts gives
$$\displaystyle
I_{c}(\tau;s)\vert_{s=ca}
\sim e^{\tau ca} \int_0^T\left(v_x(a,t)u(a,t)-u_x(a,t)v(a,t)\right)dt
\tag {1.7}
$$
modulo exponentially decaying as $\tau\longrightarrow\infty$.
Since the set $]a,\,\infty[\,\times\,]0,\,T[$ and $s=ca$
correspond to $D$ and $s=h_D(\omega)$, respectively and the right
hand side of (1.7) corresponds to the integral (1.3) in the
elliptic case as explained above, from the past experience we expected that this right hand
side decays truly algebraically as $\tau\longrightarrow\infty$. If
it is true, then we automatically obtain the formula
$$\displaystyle
\lim_{\tau\longrightarrow\infty}
\frac{\displaystyle\log\vert I_{c}(\tau;0)\vert}
{\tau}
=-ca
$$
which corresponds to the formula (1.2). However, the fact is
different from the expected value.  See Figure 1 below.

\begin{figure}
\begin{center}
\epsfxsize=9cm
\epsfysize=9cm
\epsfbox{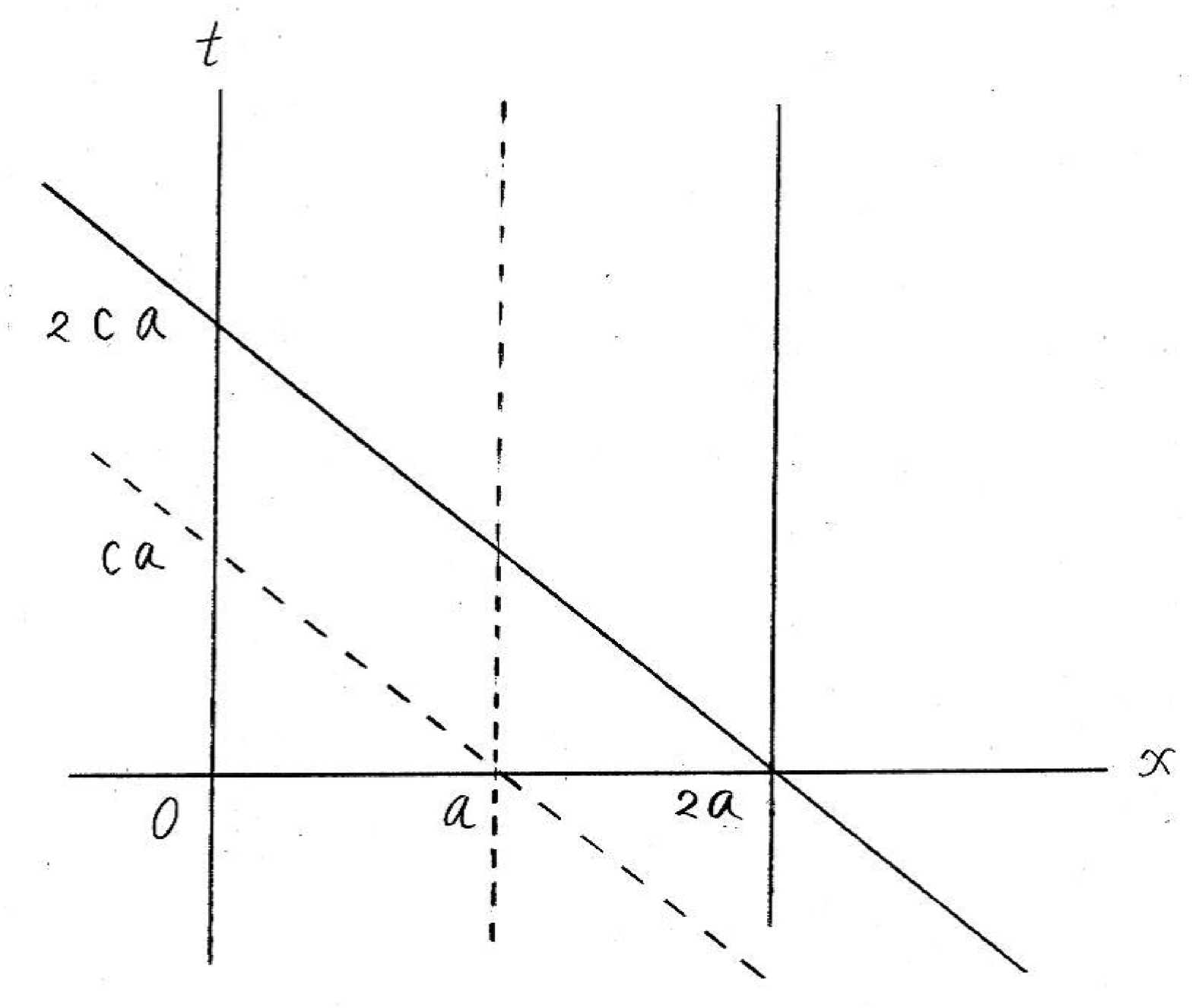}
\caption{(1.6) extracts the line $cx+t=2ca$ not $cx+t=ca$.}
\end{center}
\end{figure}

\newpage

\noindent One heuristic, however, technical explanation for this {\it
phenomenon} is the following.
One knows that $u(x,t)$ can be extended to the domain $]a,\,2a[\,\times\,]0,\,T[$
by the reflection $x$ $\longrightarrow$ $2a-x$ at $x=a$ as a solution of the heat equation.
Then integration by parts gives
$$\displaystyle
e^{\tau ca} \int_0^T\left(v_x(a,t)u(a,t)-u_x(a,t)v(a,t)\right)dt
\sim
e^{\tau ca} \int_0^T\left(v_x(2a,t)u(0,t)+u_x(0,t)v(2a,t)\right)dt.
$$
However, the function $e^{\tau ca}v(x,t)$ is exponentially
decaying as $\tau\longrightarrow\infty$ in a neighbourhood of the set
$\{2a\}\,\times\,[0,\,T]$.
Therefore one concludes the exponential
decaying of the function $I_{c}(\tau;s)\vert_{s=ca}$. This
argument {\it suggests} the existence of the phase function of $I_{c}(\tau;s)\vert_{s=ca}$ with a
negative real part.  Theorem 1.1 clarifies that it is essentially $-\tau ca$.

So the next question about (1.6) is: what is the value $2ca$? One
{\it interpretation} for this is: it gives us the {\it travel
time} of a signal with propagation speed $1/c$ that starts at the
boundary $x=0$ and the initial time $t=0$, reflects by another
boundary $x=a$ and returns to $x=0$.  See Figure 2 below.

\begin{figure}
\begin{center}
\epsfxsize=9cm 
\epsfysize=9cm 
\epsfbox{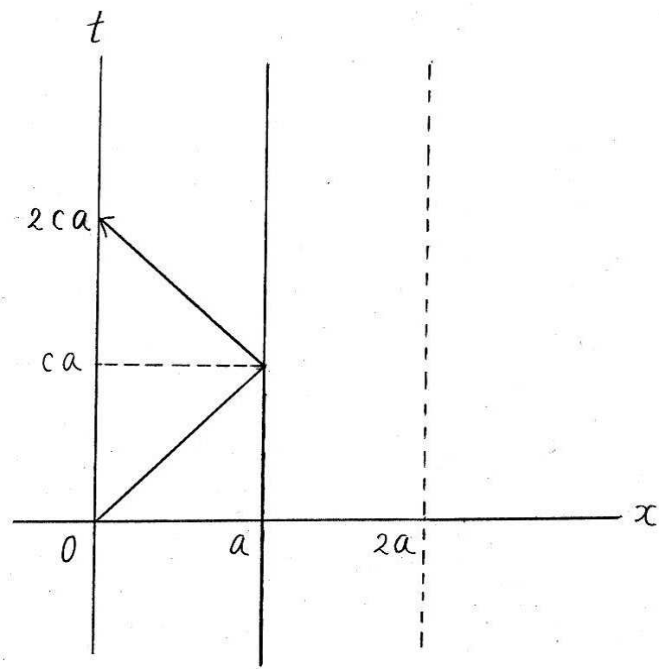}
\caption{(1.6) extracts the travel time $2ca$ of the signal with propagation speed
$1/c$ started at the boundary $x=0$.}
\end{center}
\end{figure}

\noindent
This is a quite attractive interpretation.  From this point of
view the restriction $T>Mc$ is quite reasonable since $Mc$ gives
an upper bound of the travel time $2ca$. The idea behind this
interpretation is the {\it belief}: solution of heat equation
$u_t=u_{xx}$ can be considered as a suitable superposition of
solutions of the wave equations $u_{tt}=\frac{1}{c^2}u_{xx}$ in an
appropriate sense. This belief is coming from a well known fact
for the solution of the initial value problem for the heat
equation. Given initial temperature $f=f(x), x\in\Bbb R^n$ with
compact support, let $u=u(x,t)$ be a solution of the wave equation
$u_{tt}=\triangle u$ for $x\in\Bbb R^n$, $t>0$ with initial values
$u(x,0)=f(x)$ and $u_t(x,0)=0$.  Then it is well known (e.g., see
\cite{J}) that the function
$$\displaystyle
v(x,t)=2\int_0^{\infty}\frac{e^{\displaystyle -\frac{s^2}{4t}}}{\sqrt{4\pi t}}
u(x,s)ds
$$
yields a solution of the heat equation $v_t=\triangle v$ for $x\in\Bbb R^n$, $t>0$ with
initial values $v(x,0)=f(x)$.  Then the change of variable
$s=\sqrt{t}\,\xi$ yields
$$\displaystyle
v(x,t)=
\frac{1}{\sqrt{\pi}}\,\int_0^{\infty}e^{-\xi^2/4}\,u(x,\xi\,\sqrt{t})d\xi.
$$
Replacing $t$ with $t^2$, we obtain the quite impressive formula
$$
\displaystyle
v(x,t^2)=
\frac{1}{\sqrt{\pi}}\,\int_0^{\infty}e^{-\xi^2/4}\,u(x,\xi\,t)d\xi.
$$
In one dimensional case, we have
$$
\displaystyle
v(x,t^2)=
\frac{1}{2\,\sqrt{\pi}}\,\int_0^{\infty}e^{-\xi^2/4}\,
\left(f(x+\xi\,t)+f(x-\xi\,t)\right)d\xi.
$$
The point is: the function $u(x,\xi\,t)$ of the independent
variables $x$, $t$ satisfies the wave equation with propagation
speed $\xi$, $u_{tt}=\xi^2\,\triangle u$ with initial values
$u(x,0)=f(x)$ and $u_t(x,0)=0$.  This means that $v(x,t^2)$ really
contains a signal with an arbitrary propagation speed $\xi=1/c$
and can be obtained as a {\it superposition} of signals with
several propagation speeds.  If this is true also for a solution
of (1.4), then it is reasonable to expect that the observation
data coming from the heat equation should contain some information
coming from the wave equation with arbitrary propagation speed
$\xi=1/c$.  This is another role of $c$.
In fact, in the final section we will see that a similar formula
to (1.6) is valid for the wave equation with propagation speed $1/c$.
This suggests that the indicator function for the heat equation
is a mathematical instrument that {\it picks up} a signal coming from the
corresponding wave equation with propagation speed $1/c$.

Note that, in Section 4 the reader will see that another, rather
technical restriction $T'>Mc/2$ is {\it redundant}. However, following the {\it discovery order}, we
keep the original statement since we think that the proof of
Theorem 1.1 presented in Section 3 under the restriction and the homogeneous Neumann boundary condition is still
interesting.
The proof is a time dependent approach and is based on

\noindent
$\bullet$  a representation formula of the solution of (1.4)

\noindent
$\bullet$  the explicit form of the eigenvalues and eigenfunctions
for the Laplacian with Neumann boundary condition

\noindent
$\bullet$  the asymptotic expansion of the special function
$$\displaystyle
S_1(w)=\sum_{k=0}^{\infty}(-1)^k\frac{1}{k^2+w^2},\,\,\text{Re}\,w>0
\tag {1.8}
$$
which has been established by Olver \cite{O2}.  Thus, one may think that
the proof gives an application of the special function $S_1(w)$ to an inverse problem.

\noindent
In a forthcoming paper, we will study the corresponding problems
for the heat equations with {\it variable coefficients}.

It should be pointed out that having Theorem 1.1 does not mean the
end of the inverse problem. Needless to say, the real data always
contain an error. The next problem is to consider the case when
the data are given by $u(0,t)+E_1(t)$ and $u_x(0,t)+E_2(t)$ where
the {\it size} of $E_1$ and $E_2$ are dominated by a positive
number $\delta$. Here we mean by the size of an arbitrary square
integrable function $E=E(t)$ on $]0,\,T[$ the $L^2$ norm of $E$
over $]0,\,T[$. We assume that
$$\displaystyle
\Vert(E_1,\,E_2)\Vert\equiv\Vert E_1\Vert_{L^2(0,\,T)}+\Vert E_2\Vert_{L^2(0,\,T)}\le\delta.
$$
Define the corresponding indicator function (at $s=0$) by the formula
$$\displaystyle
I_{c}(\tau;E_1, E_2)
=\int_0^T\left(-v_x(0,t)(u(0,t)+E_1(t))+(u_x(0,t)+E_2(t))v(0,t)\right)dt,\,\,\tau>c^{-2}
$$
where $v$ is the function given by (1.5).

\noindent
We assume that $\delta$ is {\it known} and {\it sufficiently small}.
The following asymptotic formula (1.10) says that if one chooses a suitable $\tau$ depending on $\delta$,
then the ratio
$$\displaystyle
\frac{\displaystyle\log\vert I_{c}(\tau;E_1,E_2)\vert}
{\tau}
$$
gives an approximate value of $-2ca$.

\proclaim{\noindent Corollary 1.2.}
Given $\delta$ and $\sigma\in\,]0,\,1[$ define
$$\displaystyle
\tau_{\sigma}(\delta)=\frac{\sigma}{T}\vert\log\delta\vert. \tag
{1.9}
$$
Then, as $\delta\longrightarrow 0$ the formula
$$\displaystyle
\sup_{\Vert(E_1,\,E_2)\Vert\le\delta}\vert\frac{\displaystyle\log\vert
I_{c}(\tau_{\sigma}(\delta);E_1,E_2)\vert}
{\tau_{\sigma}(\delta)}+2ca\vert=O\left(\frac{\vert\log\vert\log\delta\vert\vert}{\vert\log\delta\vert}\right),
\tag {1.10}
$$
is valid.
\endproclaim

\noindent Needless to say, (1.9) is one of possible choices of
$\tau$.  In practice this choice of $\tau$ will be a problem, however, it would be interested in
doing a numerical testing of the formula (1.6) as done for (1.2)
in \cite{IO}. This belongs to our future study.

A brief outline of this paper is as follows.  Theorem 1.1
is proved in Section 3.  The proof is based on an asymptotic behaviour of
an integral involving $u$ on $\{a\}\,\times\,]0,\,T'[$ in the case when
$u_x(0,t)$ is given by a constant on $]0,\,T'[$, as derived
in Section 2.
From the proof one may think that it is quite difficult to cover other boundary conditions
at $x=a$.  However, in Section 4 we present a simpler approach that is based on the transform
of $u$ into $w$ by the formula
$$\displaystyle
w(x,\tau)=\int_0^T u(x,t)e^{-z^2t}dt
$$
where
$$\displaystyle
z^2=\tau+i2c^2\tau^2\sqrt{1-\frac{1}{c^2\tau}}.
$$
We see that the asymptotic behaviour of $w(a,\tau)$ as $\tau\longrightarrow\infty$
yields a corresponding formula for the cavity with the Robin
boundary condition.  We believe that this idea will work also for the
multidimensional case.

\noindent
In the final section we give an application of
this approach to the wave equation.

\section{A key lemma.  Time dependent approach}

\noindent
The integration by parts gives
$$\displaystyle
I_c(\tau;s)
=-e^{\tau s}\int_0^Tu(a,t)v_x(a,t)dt
-e^{\tau s}\int_0^au(x,T)v(x,T)dx.
$$

\noindent
Since $v(x,t)=e^{-z^2 t}e^{xz}$, we obtain
$$\displaystyle
I_c(\tau;s)
=-ze^{\tau s}e^{az}\int_0^Tu(a,t)e^{-z^2 t}dt
-e^{\tau s}e^{-z^2T}\int_0^au(x,T)e^{xz}dx.
\tag {2.1}
$$

\noindent
Rewrite this as
$$\displaystyle
\begin{array}{c}
\displaystyle
I_c(\tau;s)
=-ze^{\tau s}e^{2az}\times e^{-az}\int_0^{T'}u(a,t)e^{-z^2 t}dt\\
\\
\displaystyle
-ze^{\tau s}e^{2az}\times e^{-az}\int_{T'}^{T}u(a,t)e^{-z^2 t}dt
-e^{\tau s}e^{-z^2T}\int_0^au(x,T)e^{xz}dx.
\end{array}
\tag {2.2}
$$

\noindent
Let $s=2ca$.  Since $T$ and $T'$ satisfy $T>2ca$
and $T'>ca$, respectively, we see that
the second and third terms of (2.2) have the estimates
$O(\tau\,e^{-\tau(T'-ca)})$ and $O(e^{-\tau(T-2ca)})$, respectively.
These yield
$$\displaystyle
I_c(\tau;s)\vert_{s=2ca}
\sim -ze^{\displaystyle -2iac\tau\sqrt{1-\frac{1}{c^2\tau}}}e^{-az}\int_0^{T'}u(a,t)e^{-z^2t}dt
\tag {2.3}
$$
modulo exponentially decaying as $\tau\longrightarrow\infty$.

\noindent Therefore it suffices to prove that
$$\displaystyle
e^{-az}\int_0^{T'}u(a,t)e^{-z^2t}dt
$$
is truly {\it algebraic decaying} as $\tau\longrightarrow\infty$.

\noindent
For the purpose we recall a known representation formula of the solution.

Let $\Omega$ be an arbitrary bounded domain of $\Bbb R^n$.
Let $u$ satisfy
$$\begin{array}{c}
\displaystyle
u_t=\triangle u\,\,\text{in}\,\Omega\times]0,\,T[,\\
\\
\displaystyle
u(x,0)=0\,\,\text{in}\,\Omega.
\end{array}
$$

\noindent
The representation formula of $u$ which is given below
involves an associated elliptic problem with the parameter $t$ and
eigenfunctions for the Laplacian with the Neumann condition.
This is taken from \cite{BC}.

\noindent Let $v=v(x,t)$ be the unique solution of the elliptic
problem depending on $t$:
$$\begin{array}{l}
\displaystyle
\triangle v=\frac{1}{\vert\Omega\vert}
\int_{\partial\Omega}\frac{\partial u}{\partial\nu}(x,t)dS(x)\,\,\text{in}\,\Omega,
\\
\\
\displaystyle
\frac{\partial v}{\partial\nu}=\frac{\partial u}{\partial\nu}(\,\cdot\,,t)\,\,\text{on}\,\partial\Omega,\\
\\
\displaystyle
\int_{\Omega}vdx=0.
\end{array}
$$

Let $\{\lambda_k\}_{k=0}^{\infty}$ and $\{\Psi_k\}_{k=0}^{\infty}$
be the all of eigenvalues and the corresponding complete
orthogonal system of the Laplacian in $\Omega$
with the Neumann boundary condition:
$$\begin{array}{l}
\displaystyle
\triangle \Psi_k+\lambda_k\Psi_k=0\,\,\text{in}\,\Omega,\\
\\
\displaystyle
\frac{\partial\Psi_k}{\partial\nu}=0\,\,\text{on}\,\partial\Omega,\\
\\
\displaystyle
\int_{\Omega}\vert\Psi_k(x)\vert^2 dx=1
\end{array}
$$
and
$$
\displaystyle
0=\lambda_0<\lambda_1\le\lambda_2\le\cdots\longrightarrow\infty.
$$
Note that $\Psi_k$ with $k\ge 1$ satisfies
$$\displaystyle
\int_{\Omega}\Psi_k(x)dx=0
$$
since $\Psi_0(x)=\vert\Omega\vert^{-1/2}$.

\proclaim{\noindent Lemma 2.1(\cite{BC}).}
The formula
$$\begin{array}{c}
\displaystyle
u(x,t)=v(x,t)
+\frac{1}{\vert\Omega\vert}
\int_0^t(\int_{\partial\Omega}\frac{\partial u}{\partial\nu}dS)dt\\
\\
\displaystyle
-\sum_{k=1}^{\infty}\left(
e^{-\lambda_k t}\int_{\Omega}v(x,0)\Psi_k(x)dx
+\int_0^t e^{-\lambda_k(t-s)}
(\int_{\Omega}v_t(x,s)\Psi_k(x)dx)ds\right)\Psi_k(x),
\end{array}
$$
is valid.

\endproclaim

\noindent
In the present case the domain $\Omega$ is given by $]0,\,a[$.
$v$ takes the form
$$\displaystyle
v(x,t)=-\frac{u_x(0,t)}{6a}(3x^2-6ax+2a^2).
$$

\noindent
$\Psi_k$ and $\lambda_k$ for $k\ge 1$ are given by
$$\displaystyle
\Psi_k(x)=\sqrt{\frac{2}{a}}\cos\,\frac{k\pi}{a}x,
\,\,\lambda_k=(\frac{k\pi}{a})^2.
$$

A direct computation yields
$$\displaystyle
\int_0^a v(x,0)\Psi_k(x) dx
=-u_x(0,0)\sqrt{\frac{2}{a}}\lambda_k^{-1}
$$
and
$$\displaystyle
\int_0^av_t(x,s)\Psi_k(x)dx
=-u_{xt}(0,s)\sqrt{\frac{2}{a}}\lambda_k^{-1}.
$$

Therefore from Lemma 2.1 we obtain
$$\begin{array}{c}
\displaystyle
u(x,t)
=v(x,t)
-\frac{1}{a}\int_0^tu_x(0,t)dt\\
\\
\displaystyle
+\sqrt{\frac{2}{a}}
\sum_{k=1}^{\infty}
\lambda_k^{-1}
\left(e^{-\lambda_k t}u_x(0,0)
+\int_0^t e^{-\lambda_k(t-s)}u_{xt}(0,s)ds\right)\Psi_k(x).
\end{array}
$$

\noindent
Letting $x\uparrow a$ in this formula we obtain the expression
of the boundary value of $u$ at $x=a$
$$\begin{array}{c}
\displaystyle
u(a,t)=\frac{a}{6}u_x(0,t)
-\frac{1}{a}\int_0^tu_x(0,t)dt\\
\\
\displaystyle
+\frac{2}{a}\sum_{k=1}^{\infty}(-1)^k
\lambda_k^{-1}
\left(e^{-\lambda_k t}u_x(0,0)
+\int_0^t e^{-\lambda_k(t-s)}u_{xt}(0,s)ds\right).
\end{array}
$$

\noindent
Now consider the simpler case
$$\noindent
u_x(0,t)=1,\,\,0<t<T'.
$$

\noindent
In this case we have
$$
\displaystyle
u(a,t)=\frac{a}{6}-\frac{t}{a}
+\frac{2}{a}\sum_{k=1}^{\infty}
(-1)^k\lambda_k^{-1}e^{-\lambda_kt},
\,0<t<T'.
\tag {2.4}
$$

The following result is the key for the enclosure method
and not trivial.

\proclaim{\noindent Lemma 2.2.}
$$\displaystyle
e^{-az}\int_0^{T'}u(a,t)e^{-z^2t}dt
$$
is algebraic decaying as $\tau\longrightarrow\infty$.
\endproclaim
{\it\noindent Proof.}
One sees that
$$\displaystyle
\int_0^{T'}e^{-\lambda_kt}e^{-z^2t}dt
=\frac{1}{\lambda_k+z^2}(1-e^{-(\lambda_k+z^2)T'}).
$$

\noindent
This gives
$$\displaystyle
e^{-az}\int_0^{T'}\sum_{k=1}^{\infty}(-1)^k\lambda_k^{-1}e^{-\lambda_k}e^{-z^2 t}dt
\sim e^{-az}\sum_{k=1}^{\infty}(-1)^k\frac{1}{\lambda_k(\lambda_k+z^2)}
$$
as $\tau\longrightarrow\infty$ modulo exponentially decaying.

\noindent
Moreover we see that
$$\displaystyle
e^{-az}\int_0^{T'}e^{-z^2 t}dt
\sim e^{-az}\frac{1}{z^2}
$$
and
$$
\displaystyle
e^{-az}\int_0^{T'}t e^{-z^2 t}dt
\sim e^{-az}\frac{1}{z^2}
$$
as $\tau\longrightarrow\infty$ modulo exponentially decaying.

\noindent
Using these and (2.4), we get
$$\begin{array}{c}
\displaystyle
e^{-az}\int_0^{T'}u(a,t)e^{-z^2 t}dt
\sim
e^{-az}\{\frac{a}{6z^2}-\frac{1}{az^4}
+\frac{2}{a}\sum_{k=1}^{\infty}(-1)^k\frac{1}{\lambda_k(\lambda_k+z^2)}\}\\
\\
\displaystyle
=e^{-az}\{\frac{a}{6z^2}-\frac{1}{az^4}
+\frac{2}{az^2}\sum_{k=1}^{\infty}(-1)^k\frac{1}{\lambda_k}
-\frac{2}{az^2}\sum_{k=1}^{\infty}(-1)^k\frac{1}{\lambda_k+z^2}\}\\
\\
\displaystyle
=e^{-az}\{\frac{a}{6z^2}-\frac{2a}{\pi^2 z^2}\sum_{k=1}^{\infty}(-1)^{k-1}\frac{1}{k^2}
+\frac{2}{az^2}\{\sum_{k=1}^{\infty}(-1)^{k-1}\frac{1}{\lambda_k+z^2}-\frac{1}{2z^2}\}\}.
\end{array}
$$

Here we cite the formula
$$\displaystyle
\sum_{k=1}^{\infty}(-1)^{k-1}\frac{1}{k^2}=\frac{\pi^2}{12}.
$$

\noindent
This yields
$$\displaystyle
e^{-az}\int_0^{T'}u(a,t)e^{-z^2t}dt
\sim
\frac{2e^{-az}}{az^2}\left(\sum_{k=1}^{\infty}(-1)^{k-1}\frac{1}{\lambda_k+z^2}
-\frac{1}{2z^2}\right).
\tag {2.5}
$$

Therefore the problem is the asymptotic behaviour of the function $f(w)$ defined by the formula
$$\displaystyle
f(w)=\sum_{k=1}^{\infty}(-1)^{k-1}\frac{1}{\lambda_k+w^2},\,\,\text{Re}\,w>0.
$$

\noindent
We know that Olver (see p. 302 of \cite{O2}) has already
studied the asymptotic behaviour of the function $S_1$ given by (1.8).
\noindent
Using the function one can write
$$
\begin{array}{c}
\displaystyle
f(w)=-(\frac{a}{\pi})^2 S_1(\frac{aw}{\pi})+\frac{1}{w^2}\\
\\
\displaystyle
=\frac{1}{2w^2}-(\frac{a}{\pi})^2R_1(\frac{aw}{\pi})
\end{array}
$$
where
$$\displaystyle
R_1(w)=S_1(w)-\frac{1}{2w^2}.
$$

\noindent
Then (2.5) becomes
$$\displaystyle
e^{-az}\int_0^{T'}u(a,t)e^{-z^2t}dt
\sim
-\frac{2ae^{-az}}{\pi^2z^2}R_1(-\frac{az}{\pi}).
\tag {2.6}
$$

Here we make use of the formula (p. 304, (7.04) of \cite{O2})
for an arbitrary fixed $\delta>0$ and $w\longrightarrow\infty$ in $\displaystyle\vert\text{arg}\,w\vert
\le\frac{\pi}{2}-\delta$:
$$\displaystyle
S_1(w)=\frac{1}{2w^2}+ \frac{2\pi}{\Gamma(1)w^{1/2}}
\sum_{j=0}^{\infty} \frac{K_{-1/2}((2j+1)\pi
w)}{(j+\frac{1}{2})^{-1/2}}
\tag {2.7}
$$
where $K_{-1/2}$ is the Macdonald function and has the asymptotic form
$$\displaystyle
K_{-1/2}(z) =(\frac{\pi}{2z})^{1/2}e^{-z}
\left(1+O(\frac{1}{z})\right),\,\, z\longrightarrow\infty
\,\,\text{in}\,\vert\text{arg}\,z\vert\le\frac{3\pi}{2}-\delta.
\tag {2.8}
$$
See p. 250 and p. 238 of \cite{O2}.

One implication of (2.7) and (2.8) is
the asymptotic formula
$$\displaystyle
R_1(w)=\frac{\pi e^{-\pi w}}{w}\left(1+O(\frac{1}{w})\right).
$$

\noindent
This yields
$$\displaystyle
e^{-az}\int_0^{T'}u(a,t)e^{-z^2t}dt
=\frac{2}{z^3}\left(1+O(\frac{1}{z})\right).
\tag {2.9}
$$

\noindent
This is the desired conclusion.

\noindent
$\Box$

\section{Proof of Theorem 1.1 and Corollary 1.2}

Consider the case when
$$\displaystyle
u_x(0,t)=\frac{t^m}{m!},\,0<t<T'.
$$

\noindent
We make use of the simple formulae:
$$\displaystyle
u_m(x,t)=\int_0^tu_{m-1}(x,t)dt,\,\,m=1, 2,\cdots,\,\,0<t<T'
$$
where $u_m$, $m=0,1,\cdots$ are the solutions corresponding to the condition
$u_x(0,t)=t^m/m!, 0<t<T'$.

\noindent Using this formulae and the initial condition
$u_m(x,0)=0$, we get
$$\begin{array}{c}
\displaystyle
\int_0^{T'}u_m(a,t)e^{-z^2t}dt
=-\frac{u_m(a,t)}{z^2}e^{-z^2 t}\vert_0^{T'}
+\frac{1}{z^2}\int_0^{T'}(u_{m})_t(a,t)e^{-z^2 t}dt\\
\\
\displaystyle
=-\frac{u_m(a,T')}{z^2}e^{-z^2 T'}
+\frac{1}{z^2}\int_0^{T'}u_{m-1}(a,t)e^{-z^2 t}dt
\end{array}
$$
and this thus yields the recurrence formula
$$\displaystyle
e^{-az}\int_0^{T'}u_m(a,t)e^{-z^2 t}dt
=\frac{e^{-az}}{z^2}\int_0^{T'}u_{m-1}(a,t)e^{-z^2 t}dt
+O\left(e^{-\tau(T'-ca)}\right).
$$

\noindent
This immediately yields
$$\displaystyle
e^{-az}\int_0^{T'}u_m(a,t)e^{-z^2 t}dt
\sim
\frac{e^{-az}}{(z^2)^m}\int_0^{T'}u_0(a,t)e^{-z^2 t}dt.
\tag {3.1}
$$

Now we are ready to prove Theorem 1.1.

\noindent
By the assumption, one may assume that $u_x(0,t)$ has the form
$$\displaystyle
u_x(0,t)=\sum_{m=0}^l\frac{\gamma_m}{m!}t^m,\,\,0<t<T'.
$$

\noindent Then the uniqueness of the direct problem gives
$$\displaystyle
u(x,t)=\sum_{m=0}^l\gamma_m u_m(x,t),\,\,0<t<T'.
$$

\noindent
From (3.1) we obtain
$$\begin{array}{c}
\displaystyle
e^{-az}\int_0^{T'}u(a,t)e^{-z^2 t}dt
=\sum_{m=0}^l\gamma_me^{-az}\int_0^{T'}u_m(a,t)e^{-z^2 t}dt\\
\\
\displaystyle
\sim\left(\sum_{m=0}^l\frac{\gamma_m}{(z^2)^m}\right)e^{-az}\int_0^{T'}u_0(a,t)e^{-z^2 t}dt.
\end{array}
\tag {3.2}
$$

Let $j=\min\,\{m\,\vert\,\gamma_m\not=0\}$.  The number $j$ satisfy $\gamma_j\not=0$
and $\gamma_m=0$ if $m<j$.

\noindent
A combination of (2.9) and (3.2) gives
$$\displaystyle
\lim_{\tau\longrightarrow\infty}z^{2j+3}e^{-az}
\int_0^{T'}u(a,t)e^{-z^2 t}dt
=2\gamma_j\,(\not=0)
$$
and from (2.3) we obtain the asymptotic formula
$$\displaystyle
\lim_{\tau\longrightarrow\infty}
z^{2(j+1)}e^{\displaystyle 2iac\tau\sqrt{1-\frac{1}{c^2\tau}}}I_c(\tau;s)\vert_{s=2ca}
=-2\gamma_j.
$$
From this formula one knows that: there exist positive numbers $\tau_0$ and $A$ such that,
for all $\tau\ge\tau_0$
$$\displaystyle
\frac{3A}{2}\ge\tau^{2(j+1)}e^{2ca\tau}\vert I_c(\tau;0)\vert
\tag {3.3}
$$
and
$$\displaystyle
\tau^{2(j+1)}e^{2ca\tau}\vert I_c(\tau;0)\vert
\ge\frac{A}{2}.
\tag {3.4}
$$
Now all the conclusions in Theorem 1.1 follow from (3.3) and (3.4).
In particular, we have
$$\displaystyle
\frac{\log\vert I_c(\tau;0)\vert}{\tau}=-2ca+2(j+1)\frac{\log\tau}{\tau}+O\left(\frac{1}{\tau}\right).
\tag {3.5}
$$

\noindent
$\Box$

\noindent
Finally we give a proof of Corollary 1.2.
We have
$$\displaystyle
I_c(\tau;E_1,E_2)-I_c(\tau;0)
=\int_0^T\left(-zE_1(t)+E_2(t)\right)e^{-z^2 t}dt.
$$
This yields the estimate
$$\displaystyle
\vert I_c(\tau;E_1,E_2)-I_c(\tau;0)\vert
\le\left(c\sqrt{\tau}+\frac{1}{\sqrt{2\tau}}\right)\delta.
\tag {3.6}
$$
A combination of (3.4) and (3.6) gives
$$\displaystyle
\vert I_c(\tau;E_1,E_2)\vert\ge \vert I_c(\tau;0)\vert
\left(1-\frac{2}{A}e^{2ca\tau}\tau^{2(j+1)}(c\sqrt{\tau}+\frac{1}{\sqrt{2\tau}})\delta\right)
\tag {3.7}
$$
and
$$\displaystyle
\vert I_c(\tau;E_1,E_2)\vert
\le\vert I_c(\tau;0)\vert
\left(1+
\frac{2}{A}e^{2ca\tau}\tau^{2(j+1)}(c\sqrt{\tau}+\frac{1}{\sqrt{2\tau}})\delta\right).
\tag {3.8}
$$

Now choose $\tau$ depending on $\delta$ in such a way that
$$\displaystyle
\frac{2}{A}e^{2ca\tau}\tau^{2(j+1)}(c\sqrt{\tau}+\frac{1}{\sqrt{2\tau}})\delta
\longrightarrow 0
\tag {3.9}
$$
as $\delta\longrightarrow 0$.  Since we have assumed that $T>2ca$, one possible choice of
such $\tau$ is to solve the equation
$$\displaystyle
e^{\tau T}=\delta^{-\sigma}.
$$
The solution is just the $\tau$ in Corollary 1.2.  In this case,
the estimation of the convergence rate of (3.9) is given by
$$\displaystyle
O\left(\vert\log\delta\vert^{2(j+1)+1/2}\delta^{1-\sigma}\right).
\tag {3.10}
$$
Now from (3.5) and (3.7) to (3.10) we obtain the desired formula
(1.10).

\noindent
$\Box$

\section{A generalization.  Stationary approach}

\noindent
The aim of this section is to

$\bullet$ present an alternative simpler method for the proof of Theorem 1.1

$\bullet$ give a generalization of Theorem 1.1.

\noindent
First we describe the problem.
Let $a>0$.
Let $u=u(x,t)$ be an arbitrary solution of
the problem:
$$\begin{array}{c}
\displaystyle
u_t=u_{xx}\,\,\text{in}\,]0,\,a[\times]0,\,T[,\\
\\
\displaystyle
u_x(a,t)+\rho u(a,t)=0\,\,\text{for}\,t\in\,]0,\,T[,\\
\\
\displaystyle
u(x,0)=0\,\,\text{in}\,]0,\,a[
\end{array}
\tag {4.1}
$$
here $\rho$ is an arbitrary fixed constant.

\noindent
Then the problem is: assume that both of $\rho$ and $a$ are {\it unknown}.
extract $a$ from $u(0,t)$ and $u_x(0,t)$ for $0<t<T$.

In this section the notion defined below plays the central role.

{\bf\noindent Definition 4.1.} Let $c>0$.  We say that a function
$f\in\, L^1(0,\,T)$ satisfies the {\it condition $c$} if
there exist positive constant $C_1$, $C_2$ and real numbers
$\tau_0(\ge c^{-2})$, $\mu_1$, $\mu_2$ such that, for all
$\tau\ge\tau_0$
$$\displaystyle
C_1\tau^{\mu_1}\le\vert\int_0^Tf(t)e^{-z^2t}dt\vert\le C_2\tau^{\mu_2}
$$
where
$$\displaystyle
z^2=\tau+i2c^2\tau^2\sqrt{1-\frac{1}{c^2\tau}}.
$$

\noindent It is easy to see that if $f(t)$ is given by a
polynomial of $t$ on an interval $]0,\,T'[\subset]0,\,T[$ that is
not identically zero, then for all $c>0$ $f(t)$ satisfies the
condition $c$. Moreover, it is not difficult to see that if $f(t)$
is smooth on $[0,\,T'[$ and $t=0$ is not a zero of $f(t)$ with
infinite order, then $f(t)$ satisfies the condition $c$.

{\bf\noindent Definition 4.2.}
Given $c>0$ define the {\it indicator function}
$I_{c}(\tau)$ by the formula
$$\displaystyle
I_{c}(\tau)
=\int_0^T\left(-v_x(0,t)u(0,t)+u_x(0,t)v(0,t)\right)dt,\,\,\tau>c^{-2}
$$
where $u$ satisfies (4.1) and $v$ is the function given by (1.5).

The following gives the answer to the problem mentioned above
and generalizes Theorem 1.1.

\proclaim{\noindent Theorem 4.1.} Assume that we know a positive
number $M$ such that $M\ge 2a$.  Let $c$ be an arbitrary positive
number satisfying $Mc<T$.  Let the $u_x(0,t)$ satisfy
the condition $c$.
Then the formula
$$\displaystyle
\lim_{\tau\longrightarrow\infty}
\frac{\displaystyle\log\vert I_{c}(\tau)\vert}
{\tau}
=-2ca,
\tag {4.2}
$$
is valid.

\endproclaim

{\it\noindent Proof.}
Introduce a new function $w$ by the formula
$$\displaystyle
w(x)=\int_0^T u(x,t)e^{-z^2t}dt,\,\,0<x<a.
$$
Note that, for simplicity of description we omitted indicating the dependence of $w$
on $\tau$ and $c$.
This function satisfies
$$\begin{array}{c}
\displaystyle
w^{''}-z^2w=e^{-z^2T}u(x,T)\,\,\text{in}\,]0,\,a[,\\
\\
\displaystyle
w'(a)+\rho w(a)=0
\end{array}
$$
and the indicator function has the expression
$$\displaystyle
I_c(\tau)=-zw(0)+w'(0).
$$
Then the integration by parts gives another expression
$$\displaystyle
I_c(\tau)=-(z+\rho)w(a)e^{az}
-e^{-z^2T}\int_0^au(\xi,T)d\xi.
$$
This yields
$$\displaystyle
I_c(\tau)e^{-2az}
=-(z+\rho)w(a)e^{-az}
+O\left(e^{-\tau(T-2ca)}\right).
\tag {4.3}
$$

Now we write $w(a)$ by using $w'(0)$.
It is easy to see that $w$ has the form
$$\displaystyle
w(x)=Ae^{xz}+Be^{-xz}
+\frac{e^{-z^2T}}{2z}
\left(e^{xz}\int_0^x u(\xi,T)e^{-\xi z}d\xi-e^{-xz}\int_0^xu(\xi,T)e^{\xi z}d\xi\right)
$$
where $A$ and $B$ are constants to be determined later.
Then the boundary conditions at $x=0$ and $x=a$ give
the system of equations
$$\begin{array}{l}
\displaystyle
zA-zB=w'(0),\\
\\
\displaystyle
(z+\rho)Ae^{az}-(z-\rho)Be^{-az}\\
\\
\displaystyle
=-\frac{e^{-z^2T}}{2}\left(
(1+\frac{\rho}{z})e^{az}\int_0^a u(\xi,T)e^{-\xi z}d\xi
+(1-\frac{\rho}{z})e^{-az}\int_0^a u(\xi,T)e^{\xi z}d\xi\right).
\end{array}
$$

\noindent
Solving this system of equations, we obtain
$$\begin{array}{l}
\displaystyle
Ae^{az}
=\frac{1}
{\displaystyle
z\left((1-\frac{\rho}{z})e^{-az}-(1+\frac{\rho}{z})e^{az}\right)}\\
\\
\displaystyle
\times
\left((1-\frac{\rho}{z})w'(0)
+\frac{e^{-z^2T}}{2}\{(1+\frac{\rho}{z})e^{2az}
\int_0^au(\xi,T)e^{-\xi z}d\xi
+(1-\frac{\rho}{z})\int_0^au(\xi,T)e^{\xi z}d\xi\}\right)\\
\\
\displaystyle
Be^{-az}

=\frac{1}
{\displaystyle
z\left((1-\frac{\rho}{z})e^{-az}-(1+\frac{\rho}{z})e^{az}\right)}\\
\\
\displaystyle
\times
\left((1+\frac{\rho}{z})w'(0)
+\frac{e^{-z^2T}}{2}\{(1+\frac{\rho}{z})
\int_0^au(\xi,T)e^{-\xi z}d\xi
+(1-\frac{\rho}{z})e^{-2az}
\int_0^au(\xi,T)e^{\xi z}d\xi\}\right).
\end{array}
$$

\noindent
A direct computation yields
$$\begin{array}{l}
\displaystyle
w(a)
=\frac{2w'(0)}
{\displaystyle
z\left((1-\frac{\rho}{z})e^{-az}-(1+\frac{\rho}{z})e^{az}\right)}\\
\\
\displaystyle
+\frac{e^{-z^2T}}
{\displaystyle
z\left((1-\frac{\rho}{z})e^{-az}-(1+\frac{\rho}{z})e^{az}\right)}
\left(\int_0^a u(\xi,T)e^{\xi z}d\xi
+\int_0^au(\xi,T)e^{-\xi z}d\xi\right)
\end{array}
$$
and thus this yields
$$\displaystyle
(z+\rho)w(a)e^{-az}
\sim 2(1+\frac{2\,\rho}{z}+O(\frac{1}{\tau^2}))w'(0)
\tag {4.4}
$$
modulo exponentially decaying as $\tau\longrightarrow\infty$.
From the assumption on $u_x(0,t)$, (4.3) and (4.4) one concludes
that, there exist positive constants $C_1'$, $C_2'$ and real
numbers $\tau_0'(\ge c^{-2})$, $\mu_1$ and $\mu_2$ such that for
all $\tau\ge\tau_0'$
$$\displaystyle
C'_1\,\tau^{\mu_1}\le\vert I_c(\tau)\vert e^{2ca\tau}\le C_2'\,\tau^{\mu_2}.
$$
This gives (4.2).

\noindent
$\Box$

\noindent
{\bf\noindent Remark 4.1.}
If once $a$ is known, then one can extract $\rho$ by the formula:
$$\displaystyle
\lim_{\tau\longrightarrow\infty}\frac{\displaystyle z(I_c(\tau)e^{-2az}+2w'(0))}{2\,w'(0)}=-2\rho.
$$

\section{Remark.  Application to the wave equation}

Finally we give a comment on an application to the wave equation.
Let $a>0$ and $c>0$.
Let $u=u(x,t)$ be an arbitrary solution of
the problem:
$$\begin{array}{c}
\displaystyle
u_{tt}=\frac{1}{c^2}u_{xx}\,\,\text{in}\,]0,\,a[\times]0,\,T[,\\
\\
\displaystyle
\frac{1}{c}u_x(a,t)+\rho u(a,t)=0\,\,\text{for}\,t\in\,]0,\,T[,\\
\\
\displaystyle
u(x,0)=0\,\,\text{in}\,]0,\,a[,\\
\\
\displaystyle
u_t(x,0)=0\,\,\text{in}\,]0,\,a[
\end{array}
\tag {5.1}
$$
here $\rho$ is an arbitrary fixed constant.  The quantity $1/c$ denotes the propagation speed
of the signal governed by the equation.  In contrast to the previous sections,
one can not change $c$.

\noindent
Then the problem is: assume that
both of $\rho$ and $a$ are {\it unknown}.
extract $a$ from $u(0,t)$ and $u_x(0,t)$ for $0<t<T$.

Introduce the function $w$ by the formula
$$\displaystyle
w(x)=\int_0^Tu(x,t)e^{-\tau t}dt,\,\,0<x<a.
$$
Again we omitted indicating the dependence of $w$ on $\tau$.
Then
$w$ satisfies
$$\begin{array}{c}
\displaystyle
\frac{1}{c^2}w^{''}-\tau^2w=e^{-\tau T}(u_t(x,T)+\tau u(x,T))\,\,\text{in}\,]0,\,a[,\\
\\
\displaystyle
\frac{1}{c}w'(a)+\rho w(a)=0.
\end{array}
$$
Now it is easy to see that
$$\begin{array}{l}
\displaystyle
w(a)
=-\frac{2w'(0)}
{\displaystyle
c\tau\left((1+\frac{\rho}{\tau})e^{ca\tau}-(1-\frac{\rho}{\tau})e^{-ca\tau}\right)}\\
\\
\displaystyle
-\frac{e^{-\tau T}}
{\displaystyle
\tau\left((1+\frac{\rho}{\tau})e^{ca\tau}-(1-\frac{\rho}{\tau})e^{-ca\tau}\right)}
\\
\\
\displaystyle
\times
\left(\int_0^a(u_t(\xi,T)+\tau u(\xi,T))e^{-c\xi\tau}d\xi
+\int_0^a(u_t(\xi,T)+\tau u(\xi,T))e^{c\xi\tau}d\xi\right).
\end{array}
\tag {5.2}
$$

{\bf\noindent Definition 5.1.}
Define the {\it indicator function}
$I(\tau)$ by the formula
$$\displaystyle
I(\tau)
=\int_0^T\left(-\frac{1}{c}v_x(0,t)u(0,t)+\frac{1}{c}u_x(0,t)v(0,t)\right)dt,\,\,\tau>0
$$
where $u$ satisfies (5.1) and $v$ is the function given by
$$\displaystyle
v(x,t)=e^{-\tau(cx+t)}.
$$

\noindent Integration by parts gives the expression
$$
\displaystyle
I(\tau)e^{2ca\tau}
=(\tau-\rho)w(a)e^{ca\tau}
-ce^{-\tau(T-2ca)}\int_0^a(u_t(\xi,T)+\tau u(\xi,T))e^{-c\xi\tau}d\xi.
$$

\noindent
Hereafter, using a similar arguments as done in Section 4 and (5.2), we obtain the following

\proclaim{\noindent Theorem 5.1.} Assume that we know a positive
number $M$ such that $M\ge 2a$.  Let $u_x(0,t)$ satisfy the
condition: there exist positive constants $C_1$, $C_2$ and real
numbers $\tau_0(>0)$, $\mu_1$, $\mu_2$ such that, for all
$\tau\ge\tau_0$
$$\displaystyle
C_1\tau^{\mu_1}\le\vert\int_0^Tu_x(0,t)e^{-\tau t}dt\vert\le C_2\tau^{\mu_2}.
\tag {5.3}
$$
Let $T>Mc$. Then the formula
$$\displaystyle
\lim_{\tau\longrightarrow\infty}
\frac{\displaystyle\log\vert I(\tau)\vert}
{\tau}
=-2ca,
\tag {5.4}
$$
is valid.

\endproclaim

\noindent Needless to say, if once $a$ is known, then one can
extract $\rho$ by the formula:
$$\displaystyle
\lim_{\tau\longrightarrow\infty}\frac{\displaystyle c\tau(I(\tau)e^{2ca\tau}+\frac{2}{c}w'(0))}{2\,w'(0)}=2\rho.
$$

The quantity $2ca$ coincides with the travel time of a
signal governed by the wave equation  with propagation speed $1/c$ which starts at the boundary
$x=0$ and initial time $t=0$, reflects another boundary $x=a$ and
returns to $x=0$.  The restriction $T>Mc$ is quite reasonable and
does not against the well known fact: the wave equation has
the {\it finite propagation property}.  The condition (5.3)
ensures that $u_x(0,t)$ can not be identically zero in an interval
$]0,T'[\subset\,]0,\,T[$. Therefore surely a signal occurs at the
initial time.  However, it should be emphasized that the formula
(5.4) makes use of the {\it averaged value} of the measured data
with an {\it exponential weight} over the observation time.  This
is a completely different idea from the well known approach in
nondestructive evaluation by sound wave: monitoring of the first
{\it arrival time} of the {\it echo}, one knows the travel time.

$$\quad$$

\centerline{{\bf Acknowledgement}}

This research was partially supported by Grant-in-Aid for
Scientific Research (C)(No.  18540160) of Japan  Society for
the Promotion of Science.  The author would like to thank Prof. Belishev
for remarks about the relationship between the response operator
and the our observation data for the heat equation in one space
dimension and comments on our results.

\vskip1cm
\noindent
e-mail address

ikehata@math.sci.gunma-u.ac.jp
\end{document}